\documentclass[leqno,12pt]{article}
\usepackage{latexsym}
 \usepackage{graphicx}
\usepackage{amsmath}
\usepackage{amssymb}
\usepackage{mathtools}
\usepackage{cite}
\usepackage{color}
\usepackage{subfig}
\usepackage[margin=1.3in, top=1.3in, bottom=1.3in]{geometry}
\usepackage{algorithm}
\usepackage{algorithmic}
\usepackage{hyperref}

\newcommand{\nc}{\newcommand}
\nc{\nt}{\newtheorem}
\nt{thm}{Theorem}[section]
\nt{cor}[thm]{Corollary}
\nt{prop}[thm]{Proposition}
\nt{lem}[thm]{Lemma}
\nt{defn}[thm]{Definition}
\nt{rem}[thm]{Remark}
\nt{exa}[thm]{Example}
\nt{ass}[thm]{Assumption}
\nt{alg}[thm]{Algorithm}
\nt{conj}[thm]{Conjecture}
\nt{claim}[thm]{Claim}
\nt{oracle}[thm]{Oracle}
\nc{\ip}[2]{\mbox{$\langle #1,#2 \rangle$}}
\nc{\pf}{\noindent{\bf Proof\ \ }}
\nc{\finpf}{\hfill{$\Box$}\linespace}
\nc{\linespace}{\vspace
{\baselineskip} \noindent}
\nc{\R}{{\mathbf R}}
\nc{\X}{{\mathbf X}}
\nc{\oR}{\overline{\R}}
\nc{\M}{\mathcal M}
\nc{\e}{\epsilon}

\nc{\Rn}{{\mathbf R}^n}
\nc{\inT}{\mbox{\rm int}\,}
\nc{\cl}{\mbox{\rm cl}\,}

\def\tto{\;{\lower 1pt \hbox{$\rightarrow$}}\kern -12pt
           \hbox{\raise 2.8pt \hbox{$\rightarrow$}}\;}
\newenvironment{myequation}{\setcounter{equation}{\value{thm}}
   \begin{equation}}{\addtocounter{thm}{1}\end{equation}}

\nc{\bmye}{\begin{myequation}}
\nc{\emye}{\end{myequation}}

\begin{document}
\title{
Lipschitz minimization and the Goldstein modulus
}
\author{Siyu Kong \and 
A.S. Lewis
\thanks{ORIE, Cornell University, Ithaca, NY.
\texttt{people.orie.cornell.edu/aslewis} 
\hspace{2cm} \mbox{~}
Research supported in part by National Science Foundation Grant DMS-2006990.}
}
\date{\today}
\maketitle

\begin{abstract}
Goldstein's 1977 idealized iteration for minimizing a Lipschitz objective fixes a distance --- the step size --- and relies on a certain approximate subgradient.  That ``Goldstein subgradient'' is the shortest convex combination of objective gradients at points within that distance of the current iterate.  A recent implementable Goldstein-style algorithm allows a remarkable complexity analysis
(Zhang et al.\ 2020), and a more sophisticated variant (Davis and Jiang, 2022) leverages typical objective geometry to force near-linear convergence.  To explore such methods, we introduce a new modulus, based on Goldstein subgradients, that robustly measures the slope of a Lipschitz function.  We relate near-linear convergence of Goldstein-style methods to linear growth of this modulus at minimizers.  We illustrate the idea computationally with a simple heuristic for Lipschitz minimization.  
\end{abstract}
\medskip

\medskip


\section{Introduction:  the Goldstein subgradient}
We consider a Euclidean space $\X$ with closed unit ball $B$, and a Lipschitz function 
$f \colon \X \to \R$ with Lipschitz constant $L \ge 0$.  Following  \cite{clarke}, the {\em Clarke subdifferential} of $f$ at any point \mbox{$x \in \X$}, denoted $\partial f(x)$, is the closed convex hull of all those vectors of the form 
$\lim_{r \to \infty} \nabla f(x_r)$ that arise from sequences $x_r \to x$.
The Clarke subdifferential is always nonempty, compact, and convex.  We say that $x$ is {\em Clarke-critical} when $0 \in \partial f(x)$.  For any radius $\e \ge 0$, following \cite{gold_eps_stat},  the {\em Goldstein subdifferential} of $f$ at $x$ is the set
\[
\partial_\e f(x) ~=~ \mbox{conv} \big( \partial f(x + \e B) \big),
\]
and it too is nonempty, compact, and convex.  The associated {\em Goldstein subgradient}, denoted $g_\e(x)$ is the shortest vector in $\partial_\e f(x)$.

As observed in \cite{gold_eps_stat}, whenever the Goldstein subgradient is nonzero, it satisfies a striking descent property:  for any current point $x \in \X$ and radius $\epsilon \ge 0$, the {\em Goldstein update}
\bmye \label{update}
g = g_\e(x) \qquad \mbox{and} \qquad x_{\mbox{\scriptsize new}} = x - \e \frac{g}{|g|},
\emye
ensures
\bmye \label{descent}
f(x_{\mbox{\scriptsize new}}) ~\le~ f(x) - \e|g|.
\emye
Iterating this update thus results in a simple conceptual minimization procedure, using a step of constant size $\e$ from one iterate to the next.  

Practical approximations of this procedure confront two fundamental challenges.
\begin{itemize} 
\item
How should we choose the step size $\e$?
\item
How should we approximate the Goldstein subgradient $g_\e(x)$?
\end{itemize}
With respect to the second challenge, the Goldstein subgradient is not computable in practice, except when the objective $f$ is simply structured, such as piecewise linear, or piecewise linear-quadratic.  However, approximation schemes have led to several recent advances.  Using one such approach, \cite{mit} presents a randomized algorithm for general Lipschitz minimization that, given any $\e >0$, finds a point $x$ satisfying $|g_\e(x)| \le \e$ using no more than $O(\e^{-4})$ calls to a specialized subgradient oracle.  This algorithm is a remarkable accomplishment, notwithstanding the slow rate, inspiring several follow-up studies:  \cite{ddmit,somit,jordan-lin-zampetakis,kornowski-shamir,jordan-deterministic,kong-lewis}.  All the more striking, then, is a recent breakthrough algorithm \cite{liwei} whose convergence is ``nearly linear'' (a term on which we expand later).  This subtle algorithm resolves both challenges above in an ingenious fashion.  

Extensive earlier literature on linear convergence in nonsmooth optimization includes the seminal work \cite{robinson-linear} on $\e$-subgradient descent methods, and was surveyed recently in \cite{atenas},  with a particular focus on bundle methods.  Those earlier works typically enumerate ``serious'' steps in procedures rather than all the subgradients needed to approximate the epsilon subdifferential for each serious step.  Crucially, the near-linear convergence result of \cite{liwei} simply counts all subgradients used.

Linear convergence results in nonsmooth optimization typically rely on growth conditions at minimizers.  For $\epsilon$-subgradient descent methods,  \cite{robinson-linear} uses an upper Lipschitz condition on the objective subdifferential.  For the Goldstein-style method of \cite{liwei}, the proof of near-linear convergence relies on a standard quadratic growth condition, along with several structural assumptions that are sophisticated, albeit generic for semi-algebraic objectives. 

Our aim here is to isolate a simple growth condition that might underly near-linear convergence of Goldstein-style subgradient methods.  We skirt the sophistication inherent in the algorithm of \cite{liwei} by decoupling the two challenges above.  We lay aside the challenge of approximating the Goldstein subgradient, imagining the following ideal oracle.

\begin{oracle}[Goldstein subgradient]  \label{oracle}
\mbox{} \\
{\rm
{\bf Input:}  a point $x \in \X$ and a radius $\e \ge 0$. \\
{\bf Output:}  the Goldstein subgradient $g_\e(x)$.
}
\end{oracle}

\noindent
We focus instead on the first challenge --- choosing the step size $\e$.  We identify a natural growth property of the objective $f$ at its minimizer that, in conjunction with a simple step-size strategy, ensures nearly linear convergence for Goldstein's iteration.  We furthermore verify the growth property for a representative class of nonsmooth objectives.

The conceptual optimization method we describe is simple but far from implementable. Nonetheless, we believe that the new growth condition underlying it may prove illuminating, both for variational analysts and for algorithm designers.

\section{The Goldstein modulus}
We begin by constructing a robust measure of the slope of $f$, starting with the following simple observation.

\begin{prop} \label{bounded}
For any $L$-Lipschitz objective function $f \colon \X \to \R$, any point $x \in \X$, and any radius $\e \ge 0$, the Goldstein subgradient satisfies $|g_\e(x)| \le L$.
\end{prop}

\noindent
The Goldstein update (\ref{update}) guarantees an objective decrease of $\e|g_\e(x)|$.  If $\e$ is small, then this guaranteed decrease is also small, being no larger than $L\e$.  On the other hand, the guaranteed decrease usually vanishes for large $\e$, because the ball $x+\e B$ then contains a Clarke critical point.  Choosing the step size $\e$ thus requires a compromise.

To address this compromise, the following definition introduces a modulus that robustly controls the size of Goldstein subgradients.  The growth condition in the definition is a central idea in our development.  We argue that this condition often holds in nonsmooth optimization, and illustrate its potential in explaining the convergence rate of Goldstein-style algorithms.

\begin{defn}
{\rm
For any Lipschitz function $f \colon \X \to \R$ and any point $x \in \X$, the {\em Goldsten modulus} is
the value
\[
\Gamma f(x) ~=~ \inf \{ \e \ge 0 : |g_\e(x)| \le \e \}.
\]
The Goldstein modulus {\em grows linearly} at a point $\bar x \in \X$ if there exists a constant 
$\alpha > 0$ such that 
\bmye \label{linear-growth}
\Gamma f(x) \ge \alpha|x - \bar x| \qquad \mbox{for all $x \in \X$ near $\bar x$}.
\emye
}
\end{defn}

We collect some elementary properties of the Goldstein modulus in the following result.

\begin{prop} \label{elementary}
The Goldstein modulus for an $L$-Lipschitz function $f \colon \X \to \R$ at a point $x \in \X$ satisfies the following properties.
\begin{enumerate}
\item[{\rm (i)}]
 $0 \le \Gamma f(x) \le L$.
 \item[{\rm (ii)}]
 \[
 |g_\e(x)| \quad 
 \left\{
 \begin{array}{ll}
 > \e & \big(\mbox{if}~ \e < \Gamma f(x)\big) \\
 \le \e & \big(\mbox{if}~ \e = \Gamma f(x)\big) \\
 < \e & \big(\mbox{if}~ \e > \Gamma f(x)\big)
 \end{array}
 \right.
 \]
\item[{\rm (iii)}]
 $\Gamma f(x) = 0$ if and only if $x$ is Clarke critical.  
\item[{\rm (iv)}]
If $\bar x$ is a Clarke critical point, then $\Gamma f(x) \le |x-\bar x|$.
\item[{\rm (v)}]
For any positive radius $\e < \Gamma f(x)$, the Goldstein update (\ref{update}) ensures an objective decrease larger than $\e^2$.
\end{enumerate}
\end{prop}

\pf
Property (i) follows from Proposition \ref{bounded}.  We next turning to property (ii).

For any value 
$\e \in \big[0, \Gamma f(x)\big)$, by definition we have $|g_\e(x)| > \e$.  Choose a sequence 
$(\e_r)$ decreasing strictly to $\Gamma f(x)$.  Let $m=1+\dim\X$.  By definition we have $|g_{\e_r}(x)| \le \e_r$ for each $r = 1,2,3,\ldots$, so there exist points $x_r^i \in x + \e_r B$, subgradients 
$y_r^i \in \partial f(x_r^i) \subset LB$, and weights $\lambda_r^i \in [0,1]$ for $i=1,2,\ldots,m$, satisfying
\[
\sum_{i=1}^m \lambda_r^i = 1 \qquad \mbox{and} \qquad \sum_{i=1}^m \lambda_r^i y_r^i = g_{\e_r}(x).
\]
After taking a subsequence, we can suppose that, for each $i=1,2,\ldots,m$, as $r \to \infty$ the points  $x_r^i$ converge to some vector $x^i \in x + \Gamma f(x) B$, the subgradients $y_r^i$ converge to some vector 
$y^i \in LB$, and the weights $\lambda_r^i$ converge to some weight $\lambda^i \in [0,1]$.  We deduce $\sum_i \lambda^i = 1$, and the vector $y = \sum_i \lambda^i y^i$ satisfies $|y| \le \Gamma f(x)$.  Since the Clarke subdifferential $\partial f$ has closed graph, we know $y^i \in \partial f(x^i)$ for each $i$, so $y \in \partial_{\Gamma f(x)} f(x)$.  We have thus proved
\[
|g_{\Gamma f(x)}(x)| ~\le~ \Gamma f(x).
\]
Finally, for any value $\e > \Gamma f(x)$ we have $\partial_\e f(x) \supset \partial_{\Gamma f(x)}(x)$, and hence
\[
|g_\e(x)| ~\le~ |g_{\Gamma f(x)}(x)| ~\le~ \Gamma f(x) ~<~ \e.
\]
Property (ii) follows.

From property (ii) we deduce $\Gamma f(x) = 0$ if and only if $|g_0(x)| = 0$, or equivalently $g_0(x)=0$, which amounts the point $x$ being Clarke critical.  Property (iv) follows from the observation
$0 \in \partial_{|x-\bar x|} f(x)$.  Property (v) follows from inequality (\ref{descent}).
\finpf

We return to the question of how to choose the radius $\e$ for the Goldstein update (\ref{update}).
By property (v) in Proposition \ref{elementary}, any choice of radius in the interval
\[
\Big[ \frac{1}{2} \Gamma f(x) \:,\: \Gamma f(x) \Big)
\]
ensures an objective decrease of at least $\frac{1}{4}\big(\Gamma f(x)\big)^2$.  Furthermore, we can use Oracle~\ref{oracle} to find such a radius by bisection search, as follows.

\begin{alg}[Goldstein modulus approximation] \label{approximation}
\label{bisection}
{\rm
\begin{algorithmic}
\STATE
\STATE {\bf input:} point $x \in \X$, Lipschitz constant $L > 0$
\IF{$g_0(x) = 0$}	
\RETURN{$0$}
\ENDIF
\STATE	$\e=\frac{1}{2}L$
\WHILE{$|g_\e(x)| \le \e$}
\STATE  $\e = \frac{1}{2}\e$
\ENDWHILE
\RETURN{$\e$}
\end{algorithmic}
}
\end{alg}

\begin{prop} \label{bisection-complexity}
For any $L$-Lipschitz function $f \colon \X \to \R$ and any point $x \in \X$, Algorithm \ref{approximation} returns a radius 
\bmye \label{estimate}
\e ~\in~ \Big[ \frac{1}{2} \Gamma f(x) \:,\: \Gamma f(x) \Big).
\emye
The number of oracles calls is $1$ if $x$ is Clarke critical, and otherwise is
\bmye \label{count}
2 + \big\lfloor \log L - \log\big(\Gamma f(x)\big) \big\rfloor.
\emye
\end{prop}

\pf
If $x$ is Clarke critical, then the first oracle call finds $g_0(x)=0$, and the algorithm returns the value $0$.  We therefore turn to the case when $x$ is not Clarke critical.  

After $1+m$ oracle calls, for $m=1,2,3,\ldots$, we have $\e = L2^{-m}$.  The algorithm terminates as soon as $|g_\e(x)| > \e$, or in other words, by Proposition \ref{elementary}, as soon as $\e < \Gamma f(x)$.  This latter condition is equivalent to $L2^{-m} < \Gamma f(x)$, or equivalently $m > \log L - \log\big(\Gamma f(x)\big)$.  The smallest such integer $m$ is 
\[
\bar m ~=~ 1 + \big\lfloor \log L - \log\big(\Gamma f(x)\big) \big\rfloor
\]
and formula (\ref{count}) follows. At termination, we know $|g_\e(x)| > \e$. We furthermore know 
$|g_{2\e}(x)| \le 2\e$:  if $\bar m > 1$, this follows from the condition in the while loop for the previous value of $\e$, and if $\bar m = 1$ it follows from the Proposition \ref{bounded}.  Property (\ref{estimate}) follows.
\finpf

\section{A simple Goldstein descent algorithm}
By combining the Goldstein descent iteration (\ref{update}) with the modulus approximation procedure, Algorithm \ref{approximation}, we arrive at the following simple algorithm.

\begin{alg}[Lipschitz minimization] \label{minimization}
{\rm
\begin{algorithmic}
\STATE
\STATE	{\bf input:}  initial point $x \in \X$, Lipschitz constant $L>0$, maximum iterations $\bar s$
\STATE	$s=0 					\hfill \mbox{}$	\COMMENT{counts Goldstein subgradient evaluations}
\LOOP
\STATE	$\e=\frac{1}{2}L$
\STATE	$g = g_\e(x)$
\STATE	$s=s+1$
\WHILE{$|g| \le \e$}			
\STATE  $\e = \frac{1}{2}\e		\hfill \mbox{}$	\COMMENT{bisection search}
\STATE	$g = g_\e(x)$
\STATE	$s=s+1$
\IF{$s > \bar s$}
\RETURN{$x$}
\ENDIF
\ENDWHILE
\STATE{$x = x - \e\frac{g}{|g|}$}
\ENDLOOP
\end{algorithmic}
}
\end{alg}

To study the complexity of Algorithm \ref{minimization}, we consider how the Goldstein modulus typically grows at a local minimizer.  That growth is at most linear, by Proposition \ref{elementary}(iv);  the following result assumes that it also satisfies a lower linear bound.

\begin{thm}[Convergence rate]~ \label{convergence-rate}
Consider an $L$-Lipschitz objective function $f \colon \X \to \R$ and a minimizer $\bar x \in \X$ at which the Goldstein modulus grows linearly.  Then, starting from any point $x_0 \in \X$ with sufficiently small initial gap $f(x_0) - f(\bar x)$, after $s$ subgradient calls to Oracle \ref{oracle}, the current point $x$ in Algorithm \ref{minimization} satisfies
\[
f(x) - f(\bar x) ~=~ O \Big( \frac{\log s}{s} \Big) \qquad \mbox{as}~ s \to \infty,
\]
where the implicit constant depends only the Lipschitz constant $L$, the initial gap, and the Goldstein modulus linear growth constant $\alpha$ in inequality (\ref{linear-growth}).
\end{thm}

\pf
Denote the current point $x$ after $k=0,1,2\ldots$ bisections searches by $x_k$.
The descent property~(\ref{descent}) guarantees
\[
f(x_k) - f(x_{k+1}) ~\ge~ \frac{1}{4} \big(\Gamma f(x_k)\big)^2.
\]
Define the objective gap $\gamma_k = f(x_k) - f(\bar x)$.  The Lipschitz and linear growth conditions imply $\Gamma f(x_k) \ge \frac{\alpha\gamma_k}{L}$, 
and hence
\[
\gamma_k - \gamma_{k+1} ~>~ \Big(\frac{\alpha\gamma_k}{2L}\Big)^2.
\]
With the change of variable $\hat\gamma = (\frac{\alpha}{2L})^2 \gamma$, we deduce
$\hat\gamma_k - \hat\gamma_{k+1} > \hat\gamma_k^2$.
Hence
\[
\frac{1}{\hat\gamma_{k+1}} ~>~ \frac{1}{\hat\gamma_k} + \frac{1}{1-\hat\gamma_k}
~>~ \frac{1}{\hat\gamma_k} + 1.
\]
By induction, we deduce
\[
\frac{1}{\hat\gamma_k} ~>~ \frac{1}{\hat\gamma_0} + k.
\]
Summarizing, for some constant $\beta$, depending only on the linear growth constant $\alpha$, the Lipschitz constant $L$, and the initial gap $\gamma_0$, we have bounded the nondecreasing sequence of gaps $\gamma_k$ in terms of the number of bisection searches $k$, by
\[
\gamma_k ~\le~ \frac{\beta}{k+1}.
\]

Turning to the number of oracle calls, we also have
\[
\Gamma f(x_k) ~\ge~ \alpha| x_k - \bar x | ~\ge~ \frac{\alpha}{L} \gamma_k.
\]
By Proposition \ref{bisection-complexity}, the bisection search at $x_k$ takes a number of oracle calls not exceeding
\[
2 + \log L - \log\big(\Gamma f(x_k)\big) ~\le~ 
2 + \log L - \log \Big(\frac{\alpha}{L} \gamma_k\Big) ~=~
\kappa - \log \gamma_k,
\]
for some constant $\kappa > 0$ depending only on $\alpha$ and $L$.  The total number $s_k$ of oracle calls accumulated before updating $x_k$ to $x_{k+1}$ therefore satisfies
\[
s_k ~\le~ \sum_{j = 0}^k (\kappa - \log \gamma_j) ~\le~ (k+1)(\kappa - \log\gamma_k) ~\le~ 
\frac{\beta(\kappa - \log\gamma_k)}{\gamma_k}.
\]

The strictly increasing continuous function $\phi \colon (0, 2^\kappa) \to (0,+\infty)$ defined by
\[
\phi(t) ~=~ \frac{t}{\beta(\kappa - \log t)}.
\]
has a strictly increasing continuous inverse $\phi^{-1} \colon (0,+\infty) \to (0, 2^\kappa)$.  Notice
\[
\lim_{s \to +\infty} \frac{s}{\log s} \phi^{-1}\Big(\frac{1}{s}\Big)
~=~ 
-\lim_{t \downarrow 0} \frac{t}{\phi(t)\log \phi(t)}
~=~ 
-\lim_{t \downarrow 0} \frac{\beta(\kappa - \log t)}{\log t - \log(\beta(\kappa - \log t))} ~=~ \beta.
\]
We deduce
\[
\phi^{-1}\Big( \frac{1}{s} \Big) ~=~ O \Big( \frac{\log s}{s} \Big),
\]
where the constant depends only on $\beta$.  Since we have proved 
$\phi(\gamma_k) \le  \frac{1}{s_k}$ for all $k=0,1,2,\ldots$, we deduce
$\gamma_k = O(\frac{\log s_k}{s_k})$, and the result follows.
\finpf

The idealized Goldstein descent complexity result in Theorem \ref{convergence-rate} is illuminating for two reasons.  First, it highlights how much stronger our imagined Oracle \ref{oracle} is than a standard subgradient oracle.  For comparison, in terms of the number $s$ of calls to standard oracles, the method of \cite{mit} for nonsmooth nonconvex objectives has complexity $O(s^{-1/4})$.  For convex objectives, the usual subgradient method has complexity $O(s^{-1/2})$, which improves to $O(\frac{1}{s})$ only in the strongly convex or smooth cases \cite{hazan-kale}.  

The second and more important reason that we present Theorem \ref{convergence-rate} is to highlight quite simply the impact on complexity of a single assumption: linear growth in the Goldstein modulus.  We will argue that this assumption often holds, even in the absence of strong convexity or smoothness.  For now, we simply illustrate the linear growth condition with three simple examples.

\begin{exa}
{\rm
Consider three functions, defined, for a constant $\alpha > 0$ and points $x \in \X$, by
\[
f_1(x) = \frac{1}{2}\alpha |x|^2, \qquad f_2(x) = \alpha |x|, \qquad f_3(x) = \frac{1}{4}\alpha |x|^4.
\]
Each function has a unique minimizer at the point $0$, and the functions $f_1$ and $f_2$ grow at least quadratically there, but $f_3$ does not.  The corresponding Goldstein modululi for the first two functions, namely
\[
\Gamma f_1(x) = \frac{\alpha}{1+\alpha}|x|, \qquad \Gamma f_2(x) = \min\{|x|,\alpha\},
\]
both grow linearly at $0$, but the third modulus satisfies
\[
\Gamma f_3(x) ~\le~ \alpha|x|^3,
\]
so its growth is slower than linear.
}
\end{exa}

\section{Nearly linear convergence}
We turn next from the simple sublinear rate guarantee in Theorem \ref{convergence-rate} to our main focus, which is linear convergence.  For many classical first-order algorithms, linear convergence is associated with quadratic growth, as described in the following definition.  This property is indispensable, both in this section, and in the remainder of this work.

\begin{defn} \label{quadratic}
{\rm
The function $f \colon \X \to \R$ {\em grows quadratically} at the point $\bar x$ when there exists a constant 
$\delta > 0$ such that
\[
f(x) - f(\bar x) ~\ge~ \frac{\delta}{2}|x-\bar x|^2 \qquad \mbox{for all $x \in \X$ near $\bar x$}.
\]
}
\end{defn}

\noindent
The relationship between linear convergence and quadratic growth is explored in detail in \cite{error_bound_d_lewis}.

Two distinct avenues to proving linear convergence suggest themselves.  We might focus on the iterates themselves, proving the existence of a constant $\theta \in (0,1)$ such that for all $r=0,1,2,\ldots$,
\bmye \label{points}
|x_{r+1} - \bar x| ~\le~ \theta |x_r - \bar x|.
\emye
Alternatively, we might consider instead the objective values, instead proving
\bmye \label{values}
f(x_{r+1}) - f(\bar x) ~\le~ \theta \big(f(x_r )- f(\bar x)\big) .
\emye
Instead of linear convergence, suppose that we are prepared to accept a somewhat relaxed rate.  Following \cite{liwei}, we say that $x_r \to \bar x$ {\em nearly linearly} if, for some exponent $m > 0$, ensuring an error $|x_r - \bar x|$ less than any small tolerance $\delta > 0$ requires only $r = O\big((\log \frac{1}{\delta})^m\big)$ iterations. (Linear convergence corresponds to the case $m=1$.)  We prove next that if each iteration satisfies {\em either} inequality (\ref{points}) {\em or} inequality (\ref{values}), then $x_r \to \bar x$ nearly linearly.

\begin{thm}[Near-linear convergence] \label{nearly}
Consider a locally Lipschitz function $f \colon \X \to \R$ that grows quadratically at a point 
$\bar x \in \X$, and a sequence of points $(x_r)$ in $X$ with initial point $x_0$ near $\bar x$.  For some constant $\theta \in (0,1)$, suppose that successive points in the sequence always satisfy 
$f(x_{r+1}) \le f(x_r)$ and one of the inequalities (\ref{points}) and (\ref{values}).
Then there exist constants $\alpha,\beta>0$ such that
\[
|x_r - \bar x| ~\le~ \alpha e^{-\beta\sqrt{r}}
\]
for all $r$.
\end{thm}

\pf
To simplify notation, we suppose $\bar x = 0$ and $f(0)=0$.  For sufficiently small initial points $x_0$, we know that there exist constants $\delta,L>0$ such that each point $x_r$ satisfies
\[
\frac{\delta}{2}|x_r|^2 \le f(x_r) \le L|x_r|.
\]
For any nonnegative integers $r,s$, if 
$f(x_{r+s}) > \theta f(x_r)$, then
\[
\frac{\theta}{L}f(x_r) ~<~ \frac{1}{L}f(x_{r+s}) ~\le~ |x_{r+s}| ~\le~ \theta^s|x_r| 
 ~\le~ \theta^s \sqrt{\frac{2f(x_r)}{\delta}}
\]
so
\[
\theta^s ~>~ \frac{\theta\sqrt{\delta}}{L\sqrt{2}} \sqrt{f(x_r)},
\]
and hence
\[
s < \mu - \nu\log\big(f(x_r)\big)
\]
for suitable constants $\mu,\nu > 0$.  We deduce
\[
f(x_{r+s}) \le \theta f(x_r) \qquad \mbox{for}~ s = \big\lceil \mu - \nu\log\big(f(x_r)\big) \big\rceil.
\]
Arguing inductively, we see, for $t=0,1,2,\ldots$,
\[
f(x_r) \le \theta^t f(x_0)
\]
for
\[
r ~\ge~ \sum_{u=0}^{t-1} \big\lceil \mu - \nu\log\big(\theta^uf(x_0)\big) \big\rceil
\]
and hence in particular for 
\begin{eqnarray*}
r 
&\ge& \sum_{u=0}^{t-1} \Big(1+\mu - \nu(u\log\theta + \log\big(f(x_0)\big) \Big) \\
&=&
t(1+\mu - \nu\log\big(f(x_0)\big) ~-~ \frac{\nu}{2}t(t-1)\log\theta.
\end{eqnarray*}
Thus there exists an integer $\psi > 0$ such that, for all $t=0,1,2\ldots$,
\[
f(x_r) ~\le~ f(x_{\psi t^2}) ~\le~ f(x_0) \theta^t 
\]
for any integer $r \ge \psi t^2$.  We deduce
\[
\frac{\delta}{2}|x_r|^2 ~\le~ f(x_r) ~\le~ f(x_0) \theta^{\big\lfloor\sqrt{\frac{1}{\psi}r}\big\rfloor}
\]
and hence
\[
2\log|x_r| ~\le~ \log \Big( \frac{2f(x_0)}{\delta} \Big) ~+~ \Big( \sqrt{\frac{1}{\psi}r} - 1 \Big)\log\theta,
\]
from which the result follows.
\finpf

To capture the two possible ways to control linear convergence, (\ref{points}) and (\ref{values}), in the context of the Goldstein iteration (\ref{update}), we make the following definition.

\begin{defn}
{\rm
A locally Lipschitz function $f \colon \X \to \R$ has {\em the Goldstein property} at a point 
$\bar x \in \X$ if there exist constants $\nu > \mu > 0$ and $\theta \in (0,1)$ such that, for any point $x \ne \bar x$ near $\bar x$, and any step size $\e$ satisfying
\[
\mu  \le \frac{\e }{|x-\bar x|} \le \nu, 
\]
the Goldstein subgradient $g = g_\e(x)$ is nonzero, and the point $x^+ = x - \e \frac{g}{|g|}$ satisfies
\[
\begin{array}{lrcl}
\mbox{either}	& f(x^+) - f(\bar x)	& \le & \theta \big( f(x) - f(\bar x) \big) \\
\mbox{or}		& |x^+ - \bar x| 		& \le & \theta |x - \bar x|.
\end{array}
\]
We call the interval $[\mu,\nu]$ a {\em proportionality bracket}.
}
\end{defn}

\begin{cor} \label{itempered}
Consider a locally Lipschitz function $f \colon \X \to \R$ that grows quadratically and has the Goldstein property at point $\bar x \in \X$, with proportionality bracket $[\mu,\nu]$.  Then any sequence of points generated iteratively from an initial point near $\bar x$ and updating iterates $x \ne \bar x$ according to the rule
\[
x ~\leftarrow~ x - \e \frac{g}{|g|} \qquad \mbox{for any}~ \e~ \mbox{satisfying}~ 
\mu  \le \frac{\e }{|x-\bar x|} \le \nu ~ \mbox{and} ~
g = g_\e(x)
\]
converges nearly linearly to $\bar x$ in the sense of Theorem~\ref{nearly}.
\end{cor}

Even assuming access to Oracle \ref{oracle} for Goldstein subgradients, the updating rule in Corollary~\ref{itempered} is not realistic in general because we do not know the distance between the current iterate $x$ and the minimizer $\bar x$.  Our strategy will be to estimate that distance from the Goldstein modulus, using the linear growth property.

\section{Robust growth relative to a manifold}
To verify linear growth of the Goldstein modulus for a locally Lipschitz function $f \colon \X \to \R$ at a point $\bar x \in \X$, we will rely on two further conditions, each of which concern some distinguished set $\M \subset \X$.  The first is a type of Lipschitz condition on the subdifferential $\partial f$.

\begin{defn} \label{lipschitz}
{\rm
The subdifferential $\partial f$ is {\em upper Lipschitz} relative to a set $\M \subset \X$ at a point 
$\bar x \in \M$ when there exists a constant $K>0$ such that all points $x,y$ near $\bar x$ with 
$x \in \M$ satisfy
\[
\partial f(y) ~ \subset ~ \partial f(x) + K|x-y|B.
\]
}
\end{defn}

The second property, which plays an important role in \cite{liwei}, describes how the value $f(x)$ grows as the point $x \in \X$ moves away from the set $\M$.  It assumes in particular that the set $\M$ is ${\mathcal C}^2$-smooth manifold around $\bar x$.  In that case, every point $x$ near $\bar x$ has a unique nearest point in $\M$, which we denote $P_\M(x)$.

\begin{defn} \label{aiming}
{\rm
If a set $\M \subset \X$ is a \mbox{${\mathcal C}^2$-smooth} manifold around a point $\bar x \in \M$, then a locally Lipschitz function $f \colon \X \to \R$ satisfies the {\em aiming condition} relative to 
$\M$ at $\bar x$ when there exists a constant $\mu > 0$ such that all points $x \in \X$ near $\bar x$ and subgradients $v \in \partial f(x)$ satisfy
\[
\ip{v}{x - P_\M(x)} ~\ge~ \mu d_\M(x).
\]
}
\end{defn}

\noindent
Some standard variational-analytic terminology \cite{clarke} helps illuminate the aiming condition.  We therefore pause to review this language.

Consider a locally Lipschitz function \mbox{$f \colon \X \to \R$}, and a point $x \in \X$.
The {\em Clarke directional derivative} of $f$ at $x$ in a direction $y \in \X$ is the quantity
\[
f^{\circ}(x;y) ~=~ \max_{v \in \partial f(x)} \ip{v}{y}.
\]
The function $f$ is {\em subdifferentially regular} at $x$ when the Clarke and classical directional derivatives agree in every direction:
\[
\lim_{t \downarrow 0} \frac{1}{t}\big( f(x+ty) - f(x) \big) ~=~ f^{\circ}(x;y) \qquad \mbox{for all}~ y \in \X.
\]
For example, sums of smooth and continuous convex functions are everywhere subdifferentially regular.
The {\em slope} of $f$ at $x$ is the quantity
\[
|\nabla f|(x) ~=~ \limsup_{y \to x}\frac{f(x)-f(y)}{|x-y|},
\]
unless $x$ is a local minimizer, in which case the slope is zero.  The following result is well known \cite[Proposition 8.5]{ioffe-variational}.

\begin{prop}[Slope and subgradients] \label{slope-subgradients}
At every point $x \in \X$, the slope of a locally Lipschitz function $f \colon \X \to \R$ satisfies the inequality
\[
|\nabla f|(x) ~\ge~ \min |\partial f(x)|,
\]
with equality if $f$ is subdifferentially regular at $x$.
\end{prop}

With this terminology in hand, we return to the aiming condition.  Definition~\ref{aiming} states that, at points $x$ outside the manifold $\M$ but near the point $\bar x$, the Clarke directional derivative of the function $f$ in the unit direction from $x$ towards its nearest point in $\M$ is uniformly negative:  
\[
f^{\circ} \Big( x \: ; \: \frac{1}{d_\M(x)}(P_\M(x) - x) \Big) ~\le~ -\mu.
\]

We can understand the aiming condition further by considering both the direction and the norm of subgradients $v \in \partial f(x)$.  First, the angle between $v$ and the direction from $x$ to its nearest point $P_\M(x)$ is uniformly larger than $\frac{\pi}{2}$:  this is precisely the ``aiming'' behavior from which the property derives its name.  Secondly, the subgradients $v$ are uniformly bounded away from zero:  
\bmye \label{sharpness2}
\liminf_{x \to \bar x,~ x \not\in \M} |\nabla f|(x) ~>~ 0,
\emye
In \cite{subgradient-curves}, property (\ref{sharpness2}) is called {\em identifiability} of the set $\M$ for the function $f$ at the point $\bar x$.  To summarize, if the aiming condition holds for $f$ relative to $\M$ at $\bar x \in \M$, then $\M$ is identifiable at $\bar x$ for $f$ and subgradients at nearby points outside $\M$ aim uniformly away from $\M$.

At any point $x \in \X$ and for any radius $\e \ge 0$, the slope and the Goldstein subgradient clearly satisfy
\[
|\nabla f|(x) ~\ge~ |g_\e(x)|,
\]
by Proposition \ref{slope-subgradients}.  With this in mind, we note that the aiming condition in fact implies a stronger property than identifiability, captured in the following crucial tool \cite[Lemma~4.1]{liwei}.

\begin{lem} \label{sharper}
The aiming condition, Definition \ref{aiming}, implies the existence of a constant $\gamma > 0$ such that
\[
\liminf_{x \to \bar x,~ x \not\in \M} |g_{\gamma d_\M(x)}(x)| ~>~ 0.
\]
\end{lem}

We next combine the conditions we need for linear growth of the Goldstein modulus in the following property.
 
\begin{defn} \label{robustly}
{\rm
If a set $\M \subset \X$ is a \mbox{${\mathcal C}^2$-smooth} manifold around a point \mbox{$\bar x \in \M$}, then 
a locally Lipschitz function $f \colon \X \to \R$ {\em grows robustly} relative to $\M$ at $\bar x$ when the following properties hold. 
\begin{itemize}
\item
$f$ grows quadratically at $\bar x$ (Definition \ref{quadratic}).
\item
$f$ is subdifferentially regular throughout $\M$.
\item
The restriction $f_\M \colon \M \to \R$ is ${\mathcal C}^2$-smooth.
\item
The subdifferential $\partial f$ is upper Lipschitz relative to $\M$ at $\bar x$ (Definition~\ref{lipschitz}).
\item
$f$ satisfies the aiming condition relative to $\M$ at $\bar x$ (Definition \ref{aiming}).
\end{itemize}
}
\end{defn}

\noindent
Two simple examples are worth keeping in mind.
\begin{exa}[Strongly convex functions]
{\rm
If the function $f$ is ${\mathcal C}^2$-smooth and strongly convex then it grows robustly at a minimizer relative to the whole space~$\X$.
}
\end{exa}

\begin{exa}[The norm]
{\rm
The norm grows robustly at $0$ relative to the manifold~$\{0\}$.  
}
\end{exa}

As a first step towards proving linear growth of the Goldstein modulus, we first observe an easier version:  linear growth of the slope. 

\begin{prop}[Local linear growth of slope]
Consider a set $\M \subset \X$ that is a \mbox{${\mathcal C}^2$-smooth} manifold around a point \mbox{$\bar x \in \M$}, and a locally Lipschitz function $f \colon \X \to \R$ that grows robustly relative to $\M$ at $\bar x$.  Then, at $\bar x$, the slope of $f$ grows linearly:  there exists a constant $\beta > 0$ such that
\bmye \label{growth}
|\nabla f|(x) ~\ge~ \beta|x-\bar x| \qquad \mbox{for all}~ x \in \X ~ \mbox{near}~ \bar x.
\emye
\end{prop}

\pf
The ${\mathcal C}^2$-smooth restriction $f_\M \colon \M \to \R$ grows quadratically at its local minimizer $\bar x$, so it is strongly convex on a neighborhood of $\bar x$ \cite[Theorem 11.21]{boumal2022intromanifolds}, and hence, by \cite[Lemma 11.28]{boumal2022intromanifolds}, its Riemannian gradient grows linearly at $\bar x$.  In other words, using the notation of \cite{{boumal2022intromanifolds}}, there exists a constant $\beta > 0$ such that all points $x \in \M$ near $\bar x$ satisfy
\[
\| \mbox{grad}(f_\M) \|_x  ~\ge~ \beta|x-\bar x|
\]
and consequently
\[
|\nabla f|(x) ~\ge~ |\nabla f_\M|(x) ~=~ \| \mbox{grad}(f_\M) \|_x  ~\ge~ \beta|x-\bar x|.
\]
Combined with the identifiability condition (\ref{sharpness2}), we deduce local linear growth of the slope.
\finpf 

\noindent
We can strengthen this result, as follows.  

\begin{thm} \label{linear-Goldstein}
If a locally Lipschitz function $f \colon \X \to \R$ grows robustly at a point, then its Goldstein modulus grows linearly there.
\end{thm}

\pf
To simplify notation, suppose that the point of interest is $\bar x = 0$.  We argue by contradiction.  If the result fails, then there exists a sequence of values 
$0 < \alpha_r \downarrow 0$ and a sequence of nonzero points $x_r \to 0$ in $\X$ such that 
$\Gamma f(x_r) < \alpha_r |x_r|$ for each $r=1,2,\ldots$.  Each corresponding Goldstein subgradient must then satisfy $|g_{\alpha_r|x_r|}(x_r)| < \alpha_r |x_r|$.  After taking a subsequence, we can suppose that the normalized vectors $\frac{x_r}{|x_r|}$ converge to some unit direction $u \in \X$.  

Following the notation of Lemma \ref{sharper} and Definitions \ref{lipschitz} and \ref{robustly}, suppose that the direction $u$ lies outside the tangent space $T$ to the manifold $\M$ at $0$.  As $r \to \infty$, we have 
$x_r = |x_r|u + o(|x_r|)$, 
and hence $d_\M(x_r) =  d_T(u)|x_r| + o(|x_r|)$.  For all large $r$ we deduce $\gamma d_\M(x_r) > \alpha_r|x_r|$ and hence
\[
|g_{\gamma d_\M(x_r)}(x_r)| ~\le~ |g_{\alpha_r|x_r|}(x_r)| ~<~ \alpha_r |x_r| ~\to~ 0 \qquad \mbox{as}~ r \to \infty,
\]
contradicting Lemma \ref{sharper}.

The direction $u$ must therefore lies in the tangent space $T$, so there exists a sequence of points $x'_r \in \M$ satisfying $x'_r - |x_r| u = o(|x_r|)$ and hence $x_r - x'_r = o(|x_r|)$ as $r \to \infty$.  The upper-Lipschitz property ensures
\begin{eqnarray*}
g_{\alpha_r|x_r|}(x_r) 
& \in & \partial_{\alpha_r|x_r|}f(x_r)
~=~ \mbox{conv} \big(\partial f(x_r + \alpha_r|x_r|B)\big) \\
& \subset & \partial f(x'_r) ~+~ \big(K\alpha_r|x_r| + o(|x_r|)\big)B,
\end{eqnarray*}
so there exist subgradients $y_r \in \partial f(x'_r)$ satisfying 
\[
| g_{\alpha_r|x_r|}(x_r) - y_r | ~\le~ K\alpha_r|x_r| + o(|x_r|).
\]
The linear growth property (\ref{growth}) along with the subdifferential regularity of $f$ at $x'_r$ shows 
$|y_r| \ge |\nabla f|(x'_r) \ge \beta|x'_r|$, so we deduce
\[
\alpha_r |x_r| ~>~ |g_{\alpha_r|x_r|}(x_r)| ~\ge~ 
\beta|x'_r| - K\alpha_r|x_r| + o(|x_r|) ~=~
\beta|x_r| - K\alpha_r|x_r| + o(|x_r|,
\]
which is a contradiction for large $r$.
\finpf

\section{Robust growth for max functions}
Classical nonlinear programming furnishes a central example of robust growth.  We consider {\em smooth max functions} $f \colon \X \to \R$, by which we mean functions having a representation of the form
\bmye \label{smooth-max}
f(x) ~=~ \max_{i \in I} f_i(x) \qquad (x \in \X)
\emye
for some family of continuously differentiable functions $f_i \colon \X \to \R$ indexed by $i$ in some finite index set $I$. Under reasonable conditions, we will show that smooth max functions grow robustly.  We will furthermore develop a Goldstein-style descent method that converges nearly linearly for such objectives.

At any point $x \in \X$, the subdifferential of the function (\ref{smooth-max}) is given by
\bmye \label{subdifferential}
\partial f(x) ~=~ \mbox{conv}\{\nabla f_i(x) : f_i(x) = f(x) \}
\emye
(see \cite{clarke}).
The standard second-order sufficient conditions in nonlinear programming motivate the following definition.

\begin{defn} \label{C2-max}
{\rm
A function $f \colon \X \to \R$ is a {\em strong ${\mathcal C}^2$ max function} at a point $\bar x \in \X$  if the following conditions hold.
\begin{itemize}
\item
The point $\bar x$ is Clarke critical and {\em nondegenerate\/}:   
$0 \in \mbox{ri}\big(\partial f(\bar x)\big)$.
\item
The function $f$ grows quadratically at $\bar x$ (Definition \ref{quadratic}).
\item
For some finite set $I$ and ${\mathcal C}^2$-smooth functions $f_i \colon \X \to \R$ (for $i\in I$), 
\[
f(x) ~=~ \max_{i \in I} f_i(x) \qquad \mbox{for all $x$ near $\bar x$},
\]
and furthermore the values $f_i(\bar x)$ (for $i \in I$) are all equal, and the gradients 
$\nabla f_i(\bar x)$ (for $i \in I$) are affinely independent.
\end{itemize}
In that case, we call the set of those points $x \in \X$ where the values $f_i(x)$ (for $i \in I$) are all equal the {\em active manifold}.
}
\end{defn}

The various ingredients of Definition \ref{C2-max} correspond exactly to standard second-order sufficient conditions for the nonlinear program
\[
\inf_{x \in \X,~ t \in \R} \{ t : f_i(x) \le t~ \mbox{for all}~ i \in I \},
\]
as presented in \cite{nocedal_wright}, for example.  Denote the active manifold by $\M$.  The classical first-order necessary optimality condition requires the existence of a Lagrange multiplier vector $\lambda \in \R^I_+$ satisfying $\sum_i \lambda_i = 1$ and $\sum_i \lambda_i \nabla f_i(\bar x) = 0$:  exactly the condition that $\bar x$ is Clarke critical.  The affine independence condition in the definition is just the usual linear independence constraint qualification, which ensures that $\lambda$ is unique, and furthermore that the set $\M$ is a ${\mathcal C}^2$-smooth manifold around $\bar x$, with tangent space
\[
T ~=~ \{\nabla f_i(\bar x) - \nabla f_j(\bar x) : i \ne j \}^{\perp}.
\]
Nondegeneracy reduces to the condition $\lambda_i > 0$ for all $i \in I$, which is the classical strict complementarity condition.  The classical theory of second-order sufficient conditions shows that this condition is equivalent to positive-definiteness of the operator $\sum_i \lambda_i \nabla^2 f_i(\bar x)$ on the subspace $T$.

In Definition \ref{C2-max}, the active manifold is well defined in the following sense.  Although it depends on the functions $f_i$ involved in the representation of the function $f$, the active manifold is identifiable for $f$ at $\bar x$, in the sense of inequality (\ref{sharpness2}), and identifiable manifolds must be locally unique around $\bar x$:  see \cite{ident}.

\begin{thm}
If $f \colon \X \to \R$ is a strong ${\mathcal C}^2$ max function at a point $\bar x \in \X$, with active manifold $\M$, then $f$ grows robustly relative to $\M$ at $\bar x$, and hence its Goldstein modulus grows  linearly at $\bar x$.
\end{thm}

\pf
Robust growth is proved in \cite{liwei}, and the result then follows.
\finpf

\noindent
The analogous result holds for any function $f$ satisfying \cite[Assumption A]{liwei}.

Definition \ref{C2-max} implies in particular that the restriction $f_\M \colon \M \to \R$ is 
${\mathcal C}^2$ smooth around the point $\bar x$.  At any nearby point $x \in \M$ , we can identify the Riemannian gradient of $f_\M$ with a vector $\nabla f_\M(x)$ in the tangent space $T_\M(x)$.  This vector has the following property.

\begin{prop} \label{shortest}
If $f \colon \X \to \R$ is a strong ${\mathcal C}^2$ max function at a point $\bar x \in \X$, with active manifold $\M$, then for all points $x \in \M$ near $\bar x$, the Riemannian gradient $\nabla f_\M(x)$ is the shortest convex combination of the set of gradients $\{\nabla f_i(x) : i \in I\}$.
\end{prop}

\pf
Consider any point $x \in \M$ near $\bar x$.  For all $i \in I$ we know $f=f_i$ throughout the manifold 
$\M$, the vector $\nabla f_i(x) - \nabla f_\M(x)$ must lie in the normal space $N_\M(x)$.  Equation (\ref{subdifferential}) implies
\[
\partial f(x) ~=~ \mbox{conv}\{\nabla f_i(x) : i \in I\} ~\subset~ \nabla f_\M(x) + N_\M(x).
\]
Furthermore, a partial smoothness argument \cite[Proposition 4.3]{shanshan} shows 
$\nabla f_\M(x) \in \partial f(x)$.  Since $\nabla f_\M(x)$ is the shortest vector in the affine subspace 
$\nabla f_\M(x) + N_\M(x)$, it must also be the shortest vector in $\partial f(x)$, so the result follows.
\finpf

The following tool is useful in what follows.

\begin{prop} \label{tool}
Consider the map $\Lambda \colon \X^k \to \X$ that maps any list of $k$ vectors to its shortest convex combination.  Then, around any affinely independent list whose image lies in the relative interior of its convex hull, the map $\Lambda$ is smooth.
\end{prop}

\pf
We map any list $v= (v^1,v^2,\ldots,v^k) \in \X^k$ to the shortest convex combination of the vectors in the list.  Denote the given affinely independent list by $\bar v$.  Then the relative interior assumption ensures 
$\Lambda(\bar v) = \sum_i \bar\lambda_i \bar v^i$ for some vector $\bar\lambda > 0$ solving the optimization problem
\[
\min_{\lambda \in \R^k_+} \Big\{ \frac{1}{2}\Big|\sum_i \lambda_i  \bar v^i \Big|^2 : \sum_i \lambda_i = 1 \Big\}
\]
Since $\bar\lambda > 0$, convexity implies that $\bar\lambda$ also solves the problem
\[
\min_{\lambda \in \R^k} \Big\{ \frac{1}{2}\Big|\sum_i \lambda_i  v^i \Big|^2 : \sum_i \lambda_i = 1 \Big\}
\]
when $v=\bar v$.  The solutions of this latter problem are characterized by the linear system in 
$(\alpha,\lambda) \in \R \times \R^k$
\bmye
\left\{
\begin{array}{rcl}
\alpha ~+~ \sum_i \ip{v^i}{v^j} \lambda_i & = & 0 \qquad (j=1,2,\ldots,k) \\
\sum_i \lambda_i &=& 1.
\end{array}
\right.
\emye
When $v=\bar v$, this square system is invertible, because
\[
\begin{array}{rcl}
\alpha ~+~ \sum_i \ip{\bar v^i}{\bar v^j} \lambda_i & = & 0 \qquad (j=1,2,\ldots,k) \\
\sum_j \lambda_j &=& 0
\end{array}
\]
implies
\[
0 ~=~ \sum_i \sum_j \ip{\bar v^i}{\bar v^j} \lambda_i \lambda_j ~=~ 
 \Big| \sum_i \lambda_i \bar v^i  \Big|^2 
\]
so $\sum_i \lambda_i \bar v^i = 0$ and hence $\alpha = 0$, and furthermore, by affine independence, 
$\lambda=0$.  We deduce that the solution $(\alpha,\lambda)$ depends smoothly on the list $v$ around 
$v = \bar v$, and furthermore satisfies $\lambda > 0$, since $\bar\lambda > 0$.  Consequently we know
$\Lambda(v) = \sum_i \lambda_i v^i$ also depends smoothly on $v$.
\finpf

We end this section by proving that small Goldstein subgradients of strong ${\mathcal C}^2$ max functions must approximate Riemannian gradients on the active manifold.

\begin{thm} \label{align}
Consider a strong ${\mathcal C}^2$ max function $f$ at a local minimizer $\bar x$ with active manifold $\M$.  Then there exists a constant $\mu > 0$ such that all Goldstein subgradients at points $x \in \X$ near $\bar x$ for small radii $\e \ge 0$ satisfy
\[
|g_\e(x)| \le \mu \qquad \Rightarrow \qquad 
\left\{
\begin{array}{rcl}
x & = & P_\M(x) + O(\e) \qquad \mbox{and} \\
g_\e(x) & = & \nabla_\M f\big(P_\M(x)\big) + O(\e).
\end{array}
\right.
\]
\end{thm}

\pf
Using the terminology of Definition \ref{C2-max}, for each $i \in I$, consider the set
\[
X_i ~=~ \{ x \in \X : f_i(x) = f(x) \}.
\]
The active manifold $\M$ is just $\cap_i X_i$. The affine independence assumption and a standard metric regularity argument shows the existence of a constant $C > 0$ such that 
\[
d_\M(x) ~\le~ C\max_i d_{X_i}(x)
\]
for all points $x \in \X$ near $\bar x$.  

For each $j \in I$, denote by $\mu_j$ the distance from zero to the set
\bmye \label{jmissing}
\mbox{conv}\{\nabla f_i(\bar x) : i \ne j \},
\emye
which is strictly positive by Definition \ref{C2-max}.  Fix any constant $\mu$ in the interval $(0,\min_i \mu_i)$, and consider any point $x$ satisfying $|g_\e(x)| \le \mu$.

We first claim that, if the point $x$ is near $\bar x$ and the radius $\e \ge 0$ is small, then 
$d_{X_j}(x) \le \e$ for all $j$, and hence $d_\M(x) \le C\e$.  To see this, we argue by contradiction.  If the claim fails, then there exists a sequence of points $x_r \to \bar x$ in $\X$ and radii $
\e_r \downarrow 0$, and elements $j_r \in I$, such that $|g_{\e_r}(x_r)| \le \mu$ and 
$d_{X_{j_r}}(x_r) > \e_r$ for all $r=1,2,3,\ldots$.  After taking a subsequence, we can suppose that some element $j \in I$ satisfies $j_r = j$ for all $r$.  Using Caratheodory's theorem, for $m=\dim\X + 1$, there exist scalars $\lambda_{ik}^r \ge 0$ and points $y_{ik}^r \in x_r + \e_rB$ for all $i \ne j$ in $I$, $k=1,2,\ldots,m$, such that
\[
g_{\e_r}(x_r) ~=~ \sum_{i \ne j} \sum_{k=1}^m \lambda_{ik}^r \nabla f_i(y_{ik}^r)
\quad \mbox{and} \quad
\sum_{i \ne j} \sum_{k=1}^m \lambda_{ik}^r ~=~ 1 \quad \mbox{for all}~ r=1,2,3,\ldots.
\]
Taking another subsequence, we can suppose the existence of the limits 
$\lambda_{ik} = \lim_r \lambda_{ik}^r \in [0,1]$ for each $i \ne j$ and $k=1,2,\ldots,m$, so some limit point $\hat g$ of the Goldstein subgradients $g_{\e_r}(x_r)$  has the form
\[
\hat g ~=~ \sum_{i \ne j} \sum_{k=1}^m \lambda_{ik} \nabla f_i(\bar x),
\qquad \mbox{where} \qquad
\sum_{i \ne j} \sum_{k=1}^m \lambda_{ik} ~=~ 1.
\]
We deduce that $\hat g$ lies in the set (\ref{jmissing}), and yet $|\hat g| \le \mu$, contradicting the definition of $\mu$.  We have thus proved our claim.

Assuming that the point $x$ is near $\bar x$ and the radius $\e \ge 0$ is small, we now know 
$d_{X_i}(x) \le \e$ for all $i \in I$, so there exists a point $y_i \in (x+\e B)\cap X_i$.  There exist scalars $\lambda_{ik} \ge 0$ and points 
$y_{ik} \in (x + \e B) \cap X_i$ for $i \in I$ and $k=1,2,\ldots,m$, such that
\[
g_{\e}(x) ~=~ \sum_{i \in I} \sum_{k=1}^m \lambda_{ik} \nabla f_i(y_{ik})
\quad \mbox{and} \quad
\sum_{i \in I} \sum_{k=1}^m \lambda_{ik} ~=~ 1.
\]
For convenience, we can suppose $y_{ik} = y_i$ whenever $\lambda_{ik} = 0$.
The nonnegative scalars $\lambda_i = \sum_k \lambda_{ik}$, for $i \in I$, sum to one.  Define
\[
g_i ~=~ 
\left\{
\begin{array}{ll}
\frac{1}{\lambda_i}  \sum_{k=1}^m \lambda_{ik} \nabla f_i(y_{ik}) & (\lambda_i > 0) \\ \\
\nabla f_i(y_i) & (\lambda_i = 0).
\end{array}
\right.
\]
For each $i \in I$, we also know
\[
g_i
~\in~ \mbox{conv}\big( \nabla f_i(x + \e B) \big) ~\subset~ \nabla f_i(x) + \e LB,
\]
for a suitable Lipschitz constant $L >0$ for the gradients $\nabla f_i$ around $\bar x$.
We have
\[
g_\e(x) ~=~ \sum_{i \in I} \lambda_i g_i ~\in~ \mbox{conv}\{ g_i : i \in I \}
~\subset~ \mbox{conv}\{ \nabla f_i(y_{ik}) : i \in I,~ k=1,2,\ldots,m \}
\] 
From its definition, $g_\e(x)$ must therefore be the shortest convex combination of the vectors $\nabla f_i(y_{ik})$, so it must also be the shortest convex combination of the vectors $g_i$, for $i \in I$.  

On the other hand, by Proposition \ref{shortest}, the Riemannian gradient $\nabla f_\M\big( P_\M(x)\big)$ is the shortest convex combination of the gradients $\nabla f_i(P_\M(x)\big)$, for $i\in I$.  Furthermore, for each $i \in I$ we have
\begin{eqnarray*}
|\nabla f_i(P_\M(x)\big) - g_i| 
&\le & 
|\nabla f_i(P_\M(x)\big) - \nabla f_i(x)| + |\nabla f_i(x) - g_i| \\
&\le& Ld_\M(x) + \e L ~\le~ L(1+C)\e.
\end{eqnarray*}
By assumption, when $x = \bar x$, the Riemannian gradient $\nabla f_\M\big( P_\M(x)\big) = \nabla f_\M(\bar x)$ is zero, which lies in the relative interior of the convex hull of the corresponding gradients $\nabla f_i(P_\M(x)\big) = \nabla f_i(\bar x)$.  The result therefore now follows by Proposition~\ref{tool}.~~
\finpf

\section{Tempered growth relative to a manifold}
Returning to our general theme, given a function $f \colon \X \to \R$ and a minimizer $\bar x$, 
we are interested in the behavior of Goldstein subgradients $g_{\e}(x)$ for nearby points $x$ and radii 
$\e > 0$.    Suppose that $f$ grows robustly relative to a manifold $\M$ at $\bar x$, so the aiming condition holds. Lemma \ref{sharper} then guarantees that, in the ``small radius'' regime when the radius $\e$ is small compared with the distance to the manifold $d_\M(x)$, the Goldstein subgradient cannot be too small.  Consequently, the Goldstein update (\ref{update}) ensures a reasonable decrease in the objective value.  We now consider, by contrast, the ``large radius'' regime, where the radius is of the same order of magnitude as the distance to the minimizer $\bar x$.

Theorem \ref{align} demonstrated, for strong ${\mathcal C}^2$ max functions, how small Goldstein subgradients of the function $f$ at points $x$ must approximate Riemannian gradients of the restriction $f_\M$ at the corresponding nearest points on the manifold $\M$.  We crystallize this general behavior in the following definition.

\begin{defn} \label{tempered-growth}
{\rm
Consider a set $\M \subset \X$ that is a ${\mathcal C}^2$-smooth manifold around a point $\bar x \in \M$, and  a locally Lipschitz function $f \colon \X \to \R$ whose restriction $f|_\M$ is ${\mathcal C}^2$-smooth.  We say that $f$ has {\em tempered growth} at $\bar x$ relative to $\M$ if, given any angle 
$\theta > 0$, for any sufficiently small $\beta > 0$ and sequences of points 
$\bar x \ne x_r \to \bar x$ and radii $\e_r \le \beta|x_r - \bar x|$ with corresponding Goldstein subgradients $g_{\e_r}(x_r)$ converging to zero, the subgradients and corresponding Riemannian gradients 
$\nabla_\M f\big(P_\M(x_r)\big)$ are eventually nonzero and subtend an angle less than $\theta$.
}
\end{defn}

\begin{thm}
Strong ${\mathcal C}^2$ max functions have tempered growth.
\end{thm}

\pf
Consider a strong ${\mathcal C}^2$ max function $f$ as in Definition \ref{C2-max}.  By Theorem \ref{align}, there exist constants $C$ and $D$ such that the projections $x'_r = P_\M(x_r)$ satisfy
\begin{eqnarray*}
|x_r - x'_r| & \le & C\e_r ~\le~ \beta C|x_r - \bar x|\\
|g_{\e_r}(x_r) - \nabla_\M f(x'_r)| & \le & D\e_r ~\le~ \beta D|x_r - \bar x|
\end{eqnarray*}
for all large $r$.  Quadratic growth (Definition \ref{quadratic}) ensures 
\begin{eqnarray*}
|\nabla_\M f(x'_r)| & \ge & \delta | x'_r - \bar x|  \\
& \ge & \delta (| x_r - \bar x| - | x_r - x'_r| ) ~\ge~ \delta(1-\beta C) (| x_r - \bar x|).
\end{eqnarray*}
We deduce
\[
\frac{|g_{\e_r}(x_r) - \nabla_\M f(x'_r)|}{|\nabla_\M f(x'_r)|}
~\le~ \frac{\beta D}{\delta(1-\beta C)} ~<~ \sin\theta,
\]
providing that $\beta$ is sufficiently small.
\finpf

\begin{thm} \label{tempered}
Consider a set $\M \subset \X$ that is a \mbox{${\mathcal C}^2$-smooth} manifold around a point \mbox{$\bar x \in \M$}, and a locally Lipschitz function $f \colon \X \to \R$ that grows robustly and has tempered growth  at $\bar x$ relative to $\M$.  Then $f$ has the Goldstein property at $\bar x$, and indeed, given any constant $\gamma \in (0,1)$, the interval $[\beta\gamma,\beta]$ is a proportionality bracket for all sufficiently small $\beta > 0$.  Consequently, any sequence of points generated iteratively from an initial point near $\bar x$ and updating iterates $x \ne \bar x$ according to the rule
\[
x ~\leftarrow~ x - \e \frac{g}{|g|} \qquad \mbox{for any}~ \e~ \mbox{satisfying}~ 
\gamma  \le \frac{\e }{\beta|x-\bar x|} \le 1 ~ \mbox{and} ~
g = g_\e(x)
\]
converges nearly linearly to $\bar x$ in the sense of Theorem~\ref{nearly}.
\end{thm}

\pf
By way of contradiction, consider any small constant $\beta > 0$, and suppose the property in Definition \ref{tempered-growth} fails.  Then there exists sequences of values $0 < \mu_r \to 0$, points 
$\bar x \ne x_r \to \bar x$, radii $\e_r$ satisfying
\[
\gamma ~\le~ \frac{\epsilon_r}{\beta|x_r - \bar x|} ~\le~ 1
\]
and Goldstein subgradients $g_r = g_{\e_r}(x_r)$, such that either $g_r = 0$ or the updates 
\[
x^+_r ~=~ x_r - \e_r\frac{g_r}{|g_r|}
\]
satisfy 
\begin{eqnarray*}
f(x^+_r) & > & f(x_r) - L\mu_r|x_r-\bar x| \\
|x^+_r - \bar x|  & > & (1-\mu_r)|x_r-\bar x|.
\end{eqnarray*}
Taking a subsequence, we can suppose
\bmye \label{alpha}
\frac{\epsilon_r}{|x_r - \bar x|} ~\to~ \mbox{some}~ \alpha \in (0,\beta].
\emye
If $\beta$ is small, so is $\alpha$, in which case $g_r$ cannot be zero infinitely often, by the definition of tempered growth.  Taking a subsequence, we can therefore suppose each $g_r$ is nonzero.

By inequality (\ref{descent}) we have
\[
f(x_r^+) ~\le~ f(x_r) - \e_r|g_r|,
\]
so $L\mu_r|x_r-\bar x| > \e_r|g_r|$ and hence by property (\ref{alpha}) we deduce $g_r \to 0$.  By Lemma~\ref{sharper}, there exists a constant $\gamma > 0$ such that $\e_r > \gamma d_\M(x_r)$ for all large $r$, so each projection $x'_r = P_\M(x_r)$ satisfies $|x_r - x'_r| < \frac{1}{\gamma}{\e_r}$, and 
the definition of tempered growth now ensures that each corresponding Riemannian gradient
$g'_r = \nabla_\M f(x'_r)$ is nonzero and, with $g_r$, subtends an arbitrarily small angle: specifically, assuming the quadratic growth condition (\ref{quadratic}), then for all sufficiently small $\beta$, we can guarantee 
\[
\Big| \frac{g_r}{|g_r|} - \frac{g'_r}{|g'_r|} \Big| ~\le~ \frac{\delta}{3}.
\]
Quadratic growth guarantees
\[
\Big< \frac{x'_r - \bar x}{|x'_r - \bar x|} \:,\: \frac{g'_r}{|g'_r|} \Big> ~\ge~ \frac{2\delta}{3},
\]
so
\[
\Big< \frac{x'_r - \bar x}{|x'_r - \bar x|} \:,\: \frac{g_r}{|g_r|} \Big> ~\ge~ \frac{\delta}{3}.
\]
Now notice
\begin{eqnarray*}
(1-\mu_r)^2|x_r-\bar x|^2
&<&
|x^+_r - \bar x|^2 \\
&=&
\Big|x_r - \e_r\frac{g_r}{|g_r|} - \bar x \Big|^2 \\
&=&
|x_r - \bar x|^2 + \e_r^2 - 2\e_r\Big< x_r - \bar x \:,\: \frac{g_r}{|g_r|} \Big> \\
&\le&
|x_r - \bar x|^2 + \e_r^2 + 2\e_r|x_r-x'_r| - 2\e_r\Big< x'_r - \bar x \:,\: \frac{g_r}{|g_r|} \Big> \\
&\le&
|x_r - \bar x|^2 + \e_r^2 + \frac{2}{\gamma}\e_r^2 
- 2\e_r|x'_r - \bar x|\Big< \frac{x'_r - \bar x}{|x'_r - \bar x|} \:,\: \frac{g_r}{|g_r|} \Big> \\
&\le&
|x_r - \bar x|^2  + \Big(1 + \frac{2}{\gamma}\Big)\e_r^2 
- \frac{2\delta}{3}\e_r|x'_r - \bar x| \\
&\le&
|x_r - \bar x|^2  + \Big(1 + \frac{2}{\gamma}\Big)\e_r^2 
- \frac{2\delta}{3}\e_r(|x_r - \bar x| - |x_r - x'_r|) \\
&\le&
|x_r - \bar x|^2  + \Big(1 + \frac{2}{\gamma} + \frac{2\delta}{3\gamma}\Big)\e_r^2 
- \frac{2\delta}{3}\e_r|x_r - \bar x|.
\end{eqnarray*}
Dividing both sides by $|x_r - \bar x|^2$ and letting $r \to \infty$ shows
\[
1 ~\le~ 1 + \Big(1 + \frac{2}{\gamma} + \frac{2\delta}{3\gamma}\Big)\alpha^2 - \frac{2\delta}{3} \alpha,
\]
which is a contradiction for all sufficiently small positive $\beta$ (and hence $\alpha$).
\finpf

The following algorithm realizes the near linear convergence property in Theorem~\ref{tempered} by approximating the distance to the minimizer using the Goldstein modulus.

\begin{alg}[Minimization for Lipschitz $f$] \label{fast}
\label{minimization2}
{\rm
\begin{algorithmic}
\STATE
\STATE	{{\bf input:}  Lipschitz constant $L$, initial point $x \in \X$, multiplier $\beta > 0$}
\FOR{$\mbox{iteration} = 1,2,3,\ldots$}
\STATE	$\e=\frac{1}{2}L$
\WHILE{$|g_\e(x)| \le \e$}
\STATE  $\e = \frac{1}{2}\e$
\ENDWHILE
\STATE	$\e = \beta\e$
\STATE	$g = g_\e(x)$
\STATE{$x = x - \e\frac{g}{|g|}$}
\ENDFOR
\end{algorithmic}
}
\end{alg}

\begin{thm}
With the assumptions of Theorem \ref{tempered}, for any sufficiently small multiplier $\beta > 0$, if the initial point $x$ is sufficiently close to the minimizer $\bar x$, then Algorithm \ref{fast} converges nearly linearly to $\bar x$ in the sense of Theorem~\ref{nearly}.
\end{thm}

\pf
After each {\tt while} loop, setting $\e = \beta\e$ ensures that the radius satisfies
\[
\frac{1}{2}\Gamma f(x) ~\le~ \frac{\e}{\beta} ~< ~\Gamma f(x),
\]
by Proposition \ref{bisection-complexity}.  By Theorem \ref{linear-Goldstein}, the Goldstein modulus grows linearly:  for some constant $\alpha > 0$, we know
\[
\alpha|x - \bar x| ~\le~ \Gamma f(x) \le |x-\bar x|
\]
where the second inequality follows from the fact that $\bar x$ is Clarke critical.  We deduce
\[
\frac{\alpha}{2} ~\le~ \frac{\e}{\beta|x- \bar x|} ~<~ 1.
\]
Providing that $\beta$ is sufficiently small, the result now follows from Theorem~\ref{tempered}. 
\finpf

\section{A Goldstein-style heuristic}
In practice, we cannot implement Algorithm \ref{minimization2} to minimize a Lipschitz function $f \colon \X \to \R$, because, given a point $x \in \X$ and a radius $\e > 0$, we cannot usually compute the Goldstein subgradient $g_\e(x)$.  To explore the effectiveness of the underlying idea --- adjusting the radius $\e$ adaptively by estimating the Goldstein modulus --- we therefore resort to approximating the Goldstein subgradient $g=g_\e(x)$, using a simple, easily implementable heuristic.  

Our approach is guided by the fundamental descent property that we noted at the outset:
\[
g = g_\e(x) \ne 0 \quad \Rightarrow \quad f\Big( x - \frac{\e}{|g|} g \Big) ~\le~ f(x) - \e|g|.
\]
The following definition relaxes that property.

\begin{defn}
{\rm
Consider a locally Lipschitz function $f \colon \X \to \R$, a point $x \in \X$, and a radius $\e > 0$
An {\em approximate Goldstein subgradient} of $f$ at $x$ is a subgradient $g \in \partial_\e f(x)$ satisfying the following property:
\[
|g| \ge \e \quad \Rightarrow \quad f\Big( x - \frac{\e}{|g|} g \Big) ~<~ f(x) - \frac{\e|g|}{2}.
\]
}
\end{defn}

We can compute approximate Goldstein subgradients almost surely using a simple but ingenious randomized procedure from \cite{mit}.  Consider any vector $g$ in the Goldstein subdifferential 
$\partial_\e f(x)$ that is not an approximate Goldstein subgradient.  Then, the shortest convex convex combination $g'$ of $g$ and any subgradient of $f$ at a point uniformly distributed between $x$ and
$x - \frac{\e}{|g|} g$ is likely to be substantially shorter than $g$.  Updating $g=g'$ and repeating eventually produces an approximate Goldstein subgradient, as shown in \cite{mit}.  We describe the procedure more formally as follows.

\begin{alg}[Approximate Goldstein subgradient for Lipschitz $f$]
\label{shrink}
{\rm
\begin{algorithmic}
\STATE
\STATE	{\bf input:}  center $x \in \X$, radius $\e > 0$
\STATE	{\bf output:} approximate Goldstein subgradient $g$
\STATE	choose $g \in \partial f(x)$
\STATE	$\gamma=|g|$
\WHILE{$\gamma \ge \e$}
\STATE	$y = x - \frac{\e}{\gamma} g$	\hfill \mbox{} \COMMENT{$g$ is not small so check descent property}
\IF{$f(x)-f(y) > \frac{\e\gamma}{2}$}
\STATE	{\bf break} \hfill \mbox{} \COMMENT{$g$ is an approximate Goldstein subgradient}
\ENDIF
\STATE	sample $z \in [x,y]$ uniformly at random
\STATE	choose $h \in \partial f(z)$
\STATE	$g = \mbox{shortest vector in}~ [g,h]$
\STATE	$\gamma = |g|$
\ENDWHILE
\RETURN{$g$}
\end{algorithmic}
}
\end{alg}

Consider any point $x \in \X$ that is not Clarke critical.  By Proposition \ref{elementary}, for all small $\e > 0$, every element of the Goldstein subdifferential $\partial_\e f(x)$ has norm at least $\e$.  In particular, the approximate Goldstein subgradient produced by Algorithm \ref{shrink} has norm at least $\e$.  Consequently, starting from any initial radius $\e > 0$, if we mimic our conceptual Algorithm \ref{minimization2} by repeatedly shrinking the radius and running Algorithm \ref{shrink}, then we eventually balance the sizes of the radius and a corresponding approximate Goldstein subgradient.  The final radius is our estimate of the Goldstein modulus $\Gamma(x)$.  We describe the procedure below, including a tolerance $\bar\e > 0$ that triggers termination if we encounter a small approximate Goldstein subgradient corresponding to a small radius.

\begin{alg}[Goldstein modulus estimation for Lipschitz $f$]
\label{selection}
{\rm
\begin{algorithmic}
\STATE
\STATE	{{\bf input:}  center $x \in \X$}, initial radius $\e > 0$, tolerance $\bar\e > 0$
\STATE	{{\bf output:} Goldstein modulus estimate $\e > 0$, \\
\mbox{} \hspace{1cm} approximate Goldstein subgradient $g$ satisfying 
$|g| \ge \e$}
\REPEAT
\STATE	$\e = \frac{1}{2}\e$
\STATE	find approximate Goldstein subgradient $g \in \partial_\e f(x)$ by Algorithm \ref{shrink}
\IF{$|g|< \e < \bar\e$}
\PRINT{``$x$ is approximately stationary''}
\STATE	{\bf terminate}
\ENDIF
\UNTIL{$|g| \ge \e$}
\RETURN{radius $\e$, subgradient $g$}
\end{algorithmic}
}
\end{alg}

We now mimic the philosophy of Algorithm \ref{minimization2}, replacing its Goldstein subgradients by their approximate versions.

\begin{alg}[Minimization for Lipschitz $f$] \label{min-for-lip}
{\rm
\begin{algorithmic}
\STATE
\STATE	{{\bf input:}  Lipschitz constant $L$, initial point $x \in \X$, \\
\mbox{} \hspace{1cm} tolerance $\bar\e > 0$, multiplier $\beta > 0$, maximum iterations $n$} \\
\STATE	{{\bf output:}  approximate stationary point $x$}
\FOR{$\mbox{iteration} = 1,2,3,\ldots,n$}
\STATE	$\e=L$
\STATE	run Algorithm \ref{selection} to set $\e = \mbox{Goldstein modulus estimate}$
\STATE	$\e = 2\beta\e$
\STATE	run Algorithm \ref{selection} again: \\
$\mbox{} \quad \bullet$ shrink radius $\e$ further, if necessary \\
$\mbox{} \quad \bullet$ find approximate Goldstein subgradient $g \in \partial_\e f(x)$.
\STATE	$x = x - \e\frac{g}{|g|}$
\ENDFOR
\RETURN{$x$}
\end{algorithmic}
}
\end{alg}

\begin{figure} 
    \centering
    \includegraphics[width=0.8\textwidth]{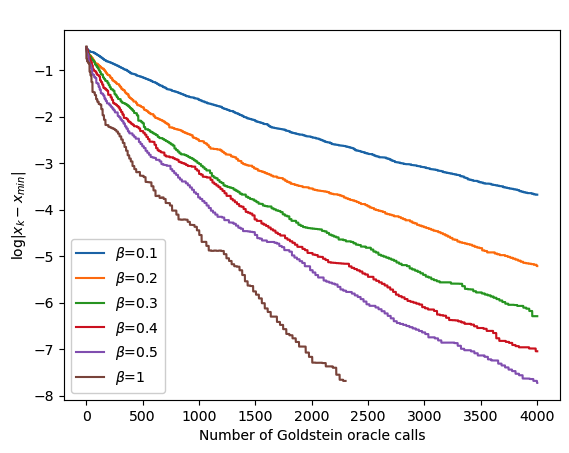}
    \caption{Algorithm \ref{min-for-lip} minimizing a maximum of five nonconvex quadratics on $\R^{10}$, using an approximate Goldstein subgradient oracle, illustrating near-linear convergence of the iterates to the minimizer.}
    \label{fig1}
\end{figure}

\begin{figure}
    \centering
    \includegraphics[width=0.8\textwidth]{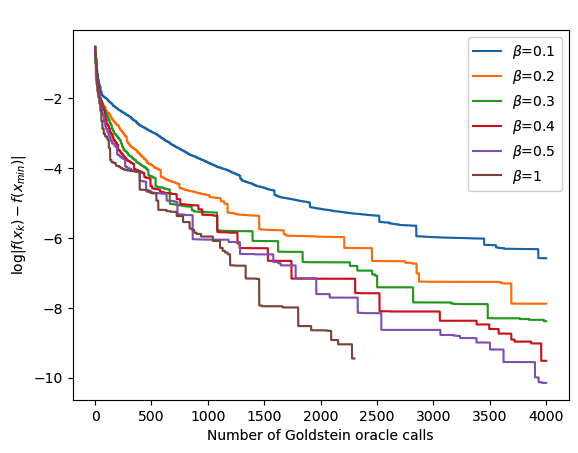}
    \caption{The example of Figure \ref{fig1}, illustrating convergence of the objective value.}
    \label{fig2}
\end{figure}

We illustrate Algorithm \ref{min-for-lip} on a simple random example.  We define a nonsmooth nonconvex function $f \colon \R^{10} \to \R$ by
\bmye \label{max-example}
f(x) ~=~ \max_{1 \le i \le 5} \{ g_i^T x + x^T H_i x \} \qquad (x \in \R^{10})
\emye
for vectors $g_i \in \R^{10}$ and $10$-by-$10$ symmetric matrices $H_i$ (for $i=1,2,3,4$) with entries uniformly distributed on the interval $[-1,1]$, and with $\sum_{i \le 5} g_i = 0$ and
$\sum_{i \le 5} H_i$ equal to the identity matrix.  This construction ensures that $f$ is almost surely a strong ${\mathcal C}^2$ max function at its global minimizer, $x_{\min} = 0$.  At any point $x \in \R^{10}$, to calculate a subgradient $h \in \partial f(x)$, we choose the index $i$ attaining the $\max$ in equation (\ref{max-example}), and set $h = g_i + 2H_i x$.  We can then run Algorithm \ref{min-for-lip} from a random initial point, for various values of the multiplier $\beta$.  We plot the progress of both the distance to the minimizer, $|x-x_{\min}|$ (in Figure \ref{fig1}) and the objective gap, $f(x) - f(x_{\min})$ (in Figure \ref{fig2}), against the number of calls to Algorithm \ref{shrink} --- our surrogate for the Goldstein subgradient oracle.

We emphasize that Algorithm \ref{min-for-lip} merely a heuristic.  We have not explored if and why the approximate Goldstein subgradients produced by Algorithm \ref{shrink} can serve as a useful substitute for the true Goldstein subgradient.  Nonetheless, the behavior of Algorithm \ref{min-for-lip}, as illustrated in Figure \ref{fig1}, is strikingly suggestive of the near-linear convergence that our theory predicts for the idealized method, Algorithm~\ref{fast}.

\begin{figure}
    \centering
    \includegraphics[width=0.8\textwidth]{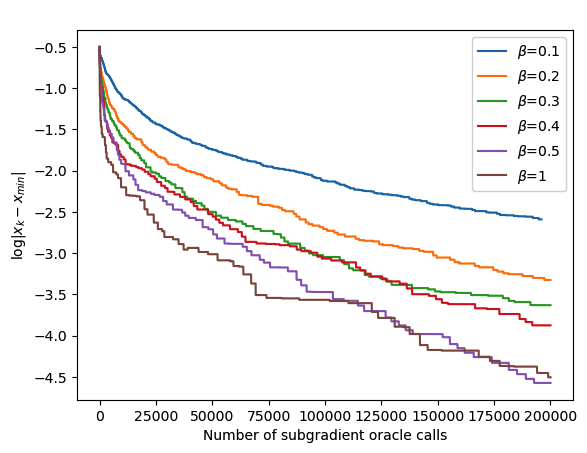}
    \caption{The example of Figure \ref{fig1}, illustrating convergence with respect to subgradient evaluations.}
    \label{fig3}
\end{figure}

Our surrogate for the Goldstein subgradient oracle, Algorithm \ref{shrink}, is written for simplicity rather than with any aim at efficiency with respect to the number of subgradients computed.  Nonetheless, for interest, Figure \ref{fig3} plots the distance to the minimizer as a function of subgradient calls.

\small
\parsep 0pt
%\bibliography{adrian}

\def\cprime{$'$} \def\cprime{$'$}

\end{document}